# Optimizing Perishable and Non-Perishable Product Assignment to packaging lines in a Sustainable Manufacturing System: An AUGMECON2VIKOR Algorithm

Reza Shahabi-Shahmiri*. Reza Tavakkoli-Moghaddam**. Zdenek Hanzalek***. Mohammad Ghasemi****. Seyed-Ali Mirnezami*****. Mohammad Rohaninejad******

*School of Industrial Engineering, College of Engineering, University of Tehran, Tehran, Iran

(e-mail: reza_shahabi011@ut.ac.ir).

**School of Industrial Engineering, College of Engineering, University of Tehran, Tehran, Iran (e-mail: tavakoli@ut.ac.ir)

*** Industrial Informatics Department, Czech Institute of Informatics Robotics and Cybernetics, Czech Technical University in Prague, Prague, Czech Republic, (e-mail: zdenek.hanzalek@cvut.cz)

**** Department of Industrial Engineering, Faculty of Engineering, Shahed University, Tehran, Iran, (e-mail: m.ghasemi8@ut.ac.ir)

***** Department of Industrial Engineering, Faculty of Engineering, Shahed University, Tehran, Iran, (e-mail: mirnezamiali@gmail.com)

****** Industrial Informatics Department, Czech Institute of Informatics Robotics and Cybernetics, Czech Technical University in Prague, Prague, Czech Republic, (e-mail: mohammad.rohani.nezhad@cvut.cz)

**Abstract**: Identifying appropriate manufacturing systems for products can be considered a pivotal manufacturing task contributing to the optimization of operational and planning activities. It has gained importance in the food industry due to the distinct constraints and considerations posed by perishable and non-perishable items in this problem. Hence, this study, proposes a new mathematical model - according to knowledge discovery as well as an assignment model to optimize manufacturing systems for perishable, non-perishable, and hybrid products tailored to meet their unique characteristics. In the presented model, three objective functions are taken into account: (1) minimizing the - production costs by assigning the products to the right set of manufacturing systems, (2) maximizing the product quality by assigning the products to the systems, and (3) minimizing the total - $CO_2$ emissions of the machines. A numerical example is utilized to evaluate the performance of AUGMECON2VIKOR compared to AUGMECON2. The results show that AUGMECON2VIKOR obtains superior Pareto solutions across all objective functions. Furthermore, the sensitivity analysis explores the positive green impacts, influencing both cost and quality.

*Keywords*: Sustainable manufacturing systems, Product features, Perishable and non-perishable products, Cost-quality-green trade-offs, AUGMECON2VIKOR algorithm.

## 1. INTRODUCTION

Over the past five decades, the manufacturing industry has thrived experiencing consistent and rapid prosperity. Manufacturers rely extensively on their ability to adapt efficiently to dynamic situations in order to survive in today's fiercely competitive production environment. Perishable products constitute a significant portion of goods produced in sectors like fresh food, pharmaceuticals, and medical services (such as blood centers). When a product spoils before customer purchase, the company may incur substantial losses (Donselaar et al., 2006). Despite the existence of various types of goods in real-world industries, in recent studies, researchers have chiefly concentrated on this type of products. At the manufacturing level, assignment problems offer a practical approach. Given the vitality of obtaining association rules considering production abilities and product features, Kuo (2019) employed data mining methods to address this problem by automatically using historical data. The objective of this reasearch was to obtain an approach to obtain an association between production abilities and product characteristics.

Production processes considering perishable items are in different sectors like the food, chemical, and pharmaceutical industries. Polotski et al. (2021) addressed this problem by extending an approach to solve the model under stochastic conditions, accounting for the unreliability of manufacturing processes, to minimize holding, backlog, and disposal cost. For investigating the production control challenges associated with perishable products, Polotski (2022) presented an approach with three steps. Initially, the hedging inventory quantity is calculated, followed by the presentation of an upper bound for the products held in stock. Ultimately, that perishable inventory bound is indicated, serving an upper bound for the hedging level. Gharbi et al. (2022) extended an optimal control approach for the manufacturing of perishable products within stochastic and dynamic production processes. Bolsi et al. (2022) explored the scheduling of perishable products and integrated workforce allocation to minimize three objective functions in a lexicographic approach. A useful heuristic, embedded in three metaheuristics was implemented. Salmasnia & Kazemi (2020) considered inventory planning, pricing, and maintenance for only perishable products using a mathematical formulation and particle swarm optimization (PSO) to maximize profit. Shah et al. (2023) presented a production inventory system in which demand depends on the average consumption of perishable goods.

To the best of our knowledge, the crucial significance of sustainable development has prompted manufacturers to consider the concepts of sustainable manufacturing. Bacher et al. (2022) proposed a sustainable flexible production system to maximize profit in which repairable defective products are reworked and products are outsourced to prevent backlogging. Castiglione et al. (2022) developed a new data formalization method for constructing a production process model, emphasizing value creation. Jafarzadeh et al. (2022) presented a mathematical formulation incorporating multiple objectives for sustainable manufacturing systems under fuzzy uncertainty to minimize total system cost, $CO_2$ emissions, and product shortages considering non-perishable products.

A review of the recent literature reveals critical gaps in addressing sustainable manufacturing systems. To fill these gaps, this article addresses the following contributions:

- Simultaneous optimization of assignment of perishable and non-perishable products to machines.
- Consideration of three machine types for perishable, non-perishable, and hybrid products.
- Investigating how assigning products to machinery affects the environment (green impacts) and improves both quality and cost
- Introduction of the AUGMECON2VIKOR algorithm to address a three-objective optimization model tailored for sustainable manufacturing systems.

## 2. PROBLEM DEFINITION

The main goal of the presented research is to introduce a novel mathematical model to improve production costs, quality, and environmental impacts by assigning the products to manufacturing systems. This assignment is carried out by discovering association relationships between manufacturing system capabilities and product features.

A food plant's packaging lines comprising three types of machines, namely perishable, non-perishable, and hybrid, are examined. Hybrid machines support both perishable and non-perishable products. In contrast with the third line, which is versatile, the first two lines are specialized and accurate. Employing the third line enhances the productivity, speed, and packaging for a specified volume of combined products, encompassing perishable and non-perishable products. It is worth noting that these combined products have a time-consuming process of separating. Additionally, there exists the possibility of exceeding the maximum capacity of the first and second lines due to a substantial amount of food products. As a result, this part of these food products can be transferred to the third line. Even though this line is applicable, given its potential capacity to accommodate a wide range of food products, it provides different quality attributes compared to the other two lines. One of the reasons necessitating the use of the third line packaging is that transferring food between the first and second lines to deal with a set of perishable and non-perishable food may result in a degradation of product quality due to temperature changes. The third line employs a Modified Atmosphere Packaging (MAP) method, reducing oxygen levels. This technique has crucial application in packaging various items such as pre-cooked foods, fish, fresh pasta, fruits, vegetables, coffee, tea, and ready-made foods. It is noteworthy that utilizing modern machinery is crucial for establishing an environmentally friendly manufacturing system because it enables advancements like better insulation and recuperation exhaust systems. These improvements help minimize energy waste and emissions, contributing to a more sustainable production process.

As illustrated in Figure 1, raw materials enter the production systems. Bluish products present perishable goods while yellowish products present non-perishable goods. Products are assigned to each type of manufacturing system with the specified capabilities. The manufacture of each one of the products within a particular manufacturing system result in a specific level of cost and quality.

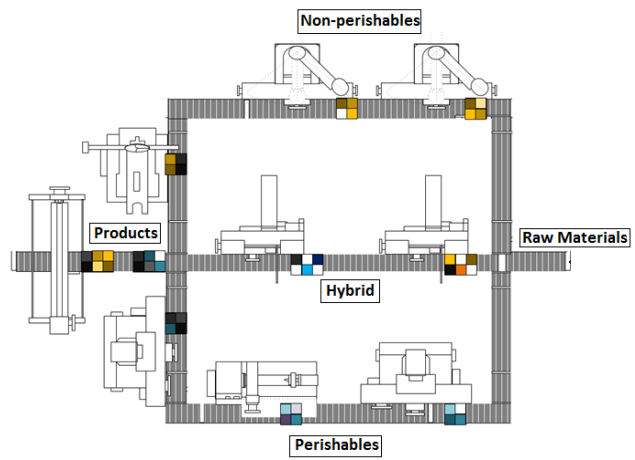

**Figure 1.** Graphical representation of the proposed model.

## 3. FORMULATION

This section presents a new mathematical programming model for assigning products to manufacturing systems. The presented model aims to find the best manufacturing systems that can fulfill all given features of the products.

*3.1 Notations*

Indices:
$n$ Non-perishable
$p$ Perishable

| | |
|---|---|
| $h$ | Perishable and non-perishable (hybrid) |
| $i$ | Product |
| $j, j'$ | Feature |
| $l$ | Manufacturing system |

Parameters:

| | |
|---|---|
| $\alpha$ | Penalty cost for exceeding the average emission |
| $\beta$ | Reward for each unit reduction below the average allowable emission |
| $AG$ | Average allowable emission |
| $G_{ijl}$ | CO2 emission amount of product $i$ with feature $j$ manufactured by system $l$ with ability to support feature $j$ |
| $Co_{ijl}$ | Cost of product $i$ with feature $j$ manufactured by system $l$ with ability to support feature $j$ |
| $q_{ijl}$ | Quality of product $i$ with feature $j$ manufactured by system $l$ with ability to support feature $j$ |
| $f_{ij}$ | 1, if Product $i$ having the $j$-th feature; otherwise 0 |
| $C_l$ | 1, if system $l$ has ability to support feature $j$; otherwise 0 |
| $d_l$ | 1, if system $l$ cannot simultaneously perform $j$th and $j'$th features; otherwise 0 |
| $Cap_l$ | Maximum capacity of product $i$ by production line $l$ |
| $R_{ijl}$ | Amount of required resources for product $i$ with feature $j$ manufactured by system $l$ with ability to support feature $j$ |
| $\gamma$ | A big number |

Binary Decision Variable:

| | |
|---|---|
| $X_{ijl}$ | 1 if product $i$ with feature $j$ is manufactured by system $l$ with ability to support fearture $j$; 0, otherwise |
| $Y1$ | 1 if CO2 emission is less than average; 0, otherwise |
| $Y2$ | 1 if CO2 emission exceeds average; 0, otherwise |

Positive Decision Variable:

| | |
|---|---|
| $EGreen$ | The amount of CO2 emission less than the average allowable emission |
| $TGreen$ | Excess CO2 emissions beyond the average |

*3.2 Equations*

$$Min\ Z_1 = \sum_{i_p}\sum_{j_p}\sum_{l_p} CO_{i_p j_p l_p} X_{i_p j_p l_p}$$
$$+ \sum_{i_n}\sum_{j_n}\sum_{l_n} CO_{i_n j_n l_n} X_{i_n j_n l_n}$$
$$+ \sum_{i_n}\sum_{j_n}\sum_{l_h} CO_{i_n j_n l_h} X_{i_n j_n l_h}$$
$$+ \sum_{i_p}\sum_{j_p}\sum_{l_h} CO_{i_p j_p l_h} X_{i_p j_p l_h}$$
$$+ TGreen.\alpha - EGreen.\beta \tag{1}$$

$$Max\ Z_2 = \sum_{i_p}\sum_{j_p}\sum_{l_p} q_{i_p j_p l_p} X_{i_p j_p l_p}$$
$$+ \sum_{i_n}\sum_{j_n}\sum_{l_n} q_{i_n j_n l_n} X_{i_n j_n l_n}$$
$$+ \sum_{i_n}\sum_{j_n}\sum_{l_h} q_{i_n j_n l_h} X_{i_n j_n l_h}$$
$$+ \sum_{i_p}\sum_{j_p}\sum_{l_h} q_{i_p j_p l_h} X_{i_p j_p l_h} \tag{2}$$

$$Min\ Z_3 = \sum_{i_p}\sum_{j_p}\sum_{l_p} G_{i_p j_p l_p} X_{i_p j_p l_p}$$
$$+ \sum_{i_n}\sum_{j_n}\sum_{l_n} G_{i_n j_n l_n} X_{i_n j_n l_n}$$
$$+ \sum_{i_n}\sum_{j_n}\sum_{l_h} G_{i_n j_n l_h} X_{i_n j_n l_h}$$
$$+ \sum_{i_p}\sum_{j_p}\sum_{l_h} G_{i_p j_p l_h} X_{i_p j_p l_h} \tag{3}$$

Subject to.

$$\sum_{l_p} C_{l_p} X_{i_p j_p l_p} + \sum_{l_h} C_{l_h} X_{i_p j_p l_h} = f_{i_p j_p} \quad ;\forall i \in i_p, \forall j \in j_p \tag{4}$$

$$\sum_{l_n} C_{l_n} X_{i_n j_n l_n} + \sum_{l_h} C_{l_h} X_{i_n j_n l_h} = f_{i_n j_n} \quad ;\forall i \in i_n, \forall j \in j_n \tag{5}$$

$$\sum_{j_p \neq j'_p} X_{i_p j_p l_p k_p} C_{l_p} \leq (1 - X_{i_p j'_p l_p m_p}).\gamma \quad \begin{array}{l};\forall i \in I_p,\\ \forall j' \in J_p,\\ \forall l \in L_p,\\ d_{l_p} = 1\end{array} \tag{6}$$

$$\sum_{j_n \neq j'_n} X_{i_n j_n l_n k_n} C_{l_n} \leq (1 - X_{i_n j'_n l_n m_n}).\gamma \quad \begin{array}{l};\forall i \in I_n,\\ \forall j' \in J_n,\\ \forall l \in L_n,\\ d_{l_n} = 1\end{array} \tag{7}$$

$$\sum_{i_n}\sum_{j_n} X_{i_n j_n l_h} + \sum_{i_p}\sum_{j_p} X_{i_p j_p l_h}$$
$$+ \sum_{i_n}\sum_{j_n} X_{i_n j_n l_h} d_{l_h}$$
$$+ \sum_{i_p}\sum_{j_p} X_{i_p j_p l_h} d_{l_h}$$
$$\leq 1 \quad ;\forall l \in L_h \tag{8}$$

$$\sum_{i_p}\sum_{j_p} X_{i_p j_p l_p} . R_{i_p j_p l_p} \leq Cap_{l_p} \quad ;\forall l \in l_p \tag{9}$$

$$\sum_{i_n}\sum_{j_n} X_{i_n j_n l_n} . R_{i_n j_n l_n} \leq Cap_{l_n} \quad ;\forall l \in l_n \tag{10}$$

$$\sum_{i_n}^{i} \sum_{j_n}^{j} X_{i_n j_n l_h} \cdot R_{i_n j_n l_h}$$
$$+ \sum_{i_p}^{i} \sum_{j_p}^{j} X_{i_p j_p l_h} \cdot R_{i_p j_p l_h} \leq Cap_{l_h} \quad ; \forall l \in L_h \quad (11)$$

$$\sum_{i_p}^{i} \sum_{j_p}^{j} \sum_{l_p}^{l} G_{i_p j_p l_p} X_{i_p j_p l_p} + \sum_{i_n}^{i} \sum_{j_n}^{j} \sum_{l_n}^{l} G_{i_n j_n l_n} X_{i_n j_n l_n}$$
$$+ \sum_{i_n}^{i} \sum_{j_n}^{j} \sum_{l_h}^{l} G_{i_n j_n l_h} X_{i_n j_n l_h} \quad (12)$$
$$+ \sum_{i_p}^{i} \sum_{j_p}^{j} \sum_{l_h}^{l} G_{i_p j_p l_h} X_{i_p j_p l_h} - AG$$
$$\leq TGreen$$

$$AG - \sum_{i_p}^{i} \sum_{j_p}^{j} \sum_{l_p}^{l} G_{i_p j_p l_p} X_{i_p j_p l_p}$$
$$+ \sum_{i_n}^{i} \sum_{j_n}^{j} \sum_{l_n}^{l} G_{i_n j_n l_n} X_{i_n j_n l_n}$$
$$+ \sum_{i_n}^{i} \sum_{j_n}^{j} \sum_{l_h}^{l} G_{i_n j_n l_h} X_{i_n j_n l_h} \quad (13)$$
$$+ \sum_{i_p}^{i} \sum_{j_p}^{j} \sum_{l_h}^{l} G_{i_p j_p l_h} X_{i_p j_p l_h} + \gamma(1 - Y1) \geq EGreen$$

$$AG - \sum_{i_p}^{i} \sum_{j_p}^{j} \sum_{l_p}^{l} G_{i_p j_p l_p} X_{i_p j_p l_p}$$
$$+ \sum_{i_n}^{i} \sum_{j_n}^{j} \sum_{l_n}^{l} G_{i_n j_n l_n} X_{i_n j_n l_n}$$
$$+ \sum_{i_n}^{i} \sum_{j_n}^{j} \sum_{l_h}^{l} G_{i_n j_n l_h} X_{i_n j_n l_h} \quad (14)$$
$$+ \sum_{i_p}^{i} \sum_{j_p}^{j} \sum_{l_h}^{l} G_{i_p j_p l_h} X_{i_p j_p l_h} \leq \gamma.Y1$$

$$\sum_{i_p}^{i} \sum_{j_p}^{j} \sum_{l_p}^{l} G_{i_p j_p l_p} X_{i_p j_p l_p} + \sum_{i_n}^{i} \sum_{j_n}^{j} \sum_{l_n}^{l} G_{i_n j_n l_n} X_{i_n j_n l_n}$$
$$+ \sum_{i_n}^{i} \sum_{j_n}^{j} \sum_{l_h}^{l} G_{i_n j_n l_h} X_{i_n j_n l_h} \quad (15)$$
$$+ \sum_{i_p}^{i} \sum_{j_p}^{j} \sum_{l_h}^{l} G_{i_p j_p l_h} X_{i_p j_p l_h} - AG$$
$$\leq \gamma.Y2$$

$$Y2 + Y1 \leq 1 \quad (16)$$
$$X_{ijl}, Y2, Y1 \in \{0,1\} \quad \forall i,j,l \quad (17)$$
$$EGreen, TGreen \geq 0$$

This model has three objective functions. The production costs are minimized in Eq. (1). Assigning each product to each machine yields a a distinct quality level, a parameter that should be maximized by the second objective in Eq. (2). Eq. (3) minimizes the total CO2 emission generated by machines for manufacturing of products. Eq. (4) states that perishable products and features can be produced with either a hybrid machine or a dedicated machine for perishable products. Eq. (5) dictates that non-perishable products and features can be produced with either a hybrid machine or a dedicated machine for non-perishable goods. Eq. (6) and (7) state that each specialized type of machine is exclusively assigned to products requiring at least one of the production features of that machine. Eq. (8) shows that each ability of a hybrid machine can support only one feature of different types of products. Eqs. (9) indicates that certain production ability of machines cannot occur simultaneously. The capacity of each machine is addressed by Eqs. (10) and (11) for perishable, non-perishable, and hybrid machines, respectively. Eqs. (12) and (13) calculate the reward and penalty costs related to the green impacts of the manufacturing system. Eqs. (14)-(16) determine whether reward and penalty costs are considered for the overall evaluation of the manufacturing system.

## 4. METHODOLOGY

### 4.1. AUGMECON2 algorithm

Marvatos and Florios (2013) introduced an augmented epsilon-constraint (AUGMECON) algorithm, incorporating a novel part $eps \times \left(\frac{S_2}{r_2} + 10^{-1}\frac{S_3}{r_3} + \cdots + 10^{-(p-2)}\frac{S_p}{r_p}\right)$ for the epsilon-constraint approach. Compared to the Epsilon constraint approach, the AUGMECON algorithm ensures the Pareto optimality of the generated solutions (Marvatos and Florios, 2013). However, growing the number of objectives leads to an escalation in processing time. Mavrotas and Florios (2013) presented an improved version of AUGMECON denoted as AUGMECON2. In this algorithm, slack/surplus variables are utilized and the bypass coefficient is considered to eliminate superfluous iterations. The reformulated objective in AUGMECON2 can be as follows:

$$Max\left(f_1(x) + eps \times \left(\frac{S_2}{r_2} + 10^{-1}\frac{S_3}{r_3} + \cdots + 10^{-(p-2)}\frac{S_p}{r_p}\right)\right) \quad (18)$$

where the parameters $r_2, r_3 \ldots, r_p$ are the ranges of the respective objective functions, $s_2, s_3 \ldots, s_p$ are the surplus variables of the respective constraints and $eps \in [10^{-6}, 10^{-3}]$.

### 4.2. AUGMECON2VIKOR approach

This algorithm has been fully described in Shahabi-Shahmiri et al. (2021). The steps of this algorithm are as follows.
Step 1: Applying VIKOR method
1.1. Construct MILP model 1 as follows:

$$Model\ 1: \begin{cases} Min\ Z_1 \\ Max\ Z_2 \\ Min\ Z_3 \\ s.t. \\ model\ constraints \end{cases}$$

1.2. Solve the model 1 and form a payoff table.
1.3. Set $f_j^*$ as
$$f_j^* = max_i(f_{ij}), f_j^- = min_i(f_{ij})$$
1.4. Consider decision maker's preferences ($w_j$) and set $S_i$ and $R_i$ as follows:

$$S_i = \sum_{j=1}^{n} w_j \, (f_j^* - f_{ij})/(f_j^* - f_j^-)$$

$$R_i = max_j[\sum_{j=1}^{n} w_j \frac{f_j^* - f_{ij}}{f_j^* - f_j^-}]$$

Step 2: Implementing augmented epsilon constraint algorithm

2.1. Construct the bi-objective mathematical model 2 as follows:

$$Model\ 2: \begin{cases} Min\ R \\ Max\ S \\ s.t. \\ model\ constraints \end{cases}$$

2.2. Solve the model 2 and form a payoff table.
2.3. Consider objective $R$ as the main objective and transfer objective $S$ to the constraints.
2.4. Calculate the range of objective $S$ ($S_{range}$).
2.5. Determine the number of grid points ($N$).
2.6. Set epsilon ($e$) as $e = \frac{S_{range}}{N}$.
2.7. Construct the model 3 as follows:

$$Model\ 3: \begin{cases} Min\ R + eps(slack/S_{range}) \\ s.t. \\ S + slack = S_{min} + p.e \\ model\ constraints \end{cases}$$

Step 3: Solving the problem
3.1. Assume that $p = 1$
3.2. Solve the model 3
3.3. If the solution is unique, add it to the pareto set and consider step 3.4; otherwise consider step 3.4.
3.4. If $p > N$ end the algorithm and report the pareto set; otherwise set $p = p + 1$ and return to step 3.2.

## 5. EXPERIMENTAL RESULTS

*5.1. Results*

A numerical example has been presented to illustrate the efficacy of the proposed model in addressing a manufacturing decision-making problem, which is crucial for corporations attempting to adapt to dynamic market and customer demands. Six products are assumed to be assigned to the production machines. There are six product features and six production abilities. Three of features and abilities are for perishable goods and the other features and abilities are for non-perishable goods. To validate the results of the AUGMECON2VIKOR algorithm and compare them to the AUGMECON2 algorithm, the presented model is coded in GAMS 24.1.2 considering five grid points. The CPU time required for AUGMECON2VIKOR to find Pareto optimal solutions was 326 seconds. Table 1 presents Pareto optimal solutions for both algorithms. Finding eight Pareto solution points, lower cost, and improved quality values in each grid point demonstrates the superiority of the AUGMECON2VIKOR algorithm. For more clarity, Fig. 2 depicts the assignment of different types of products to each machine.

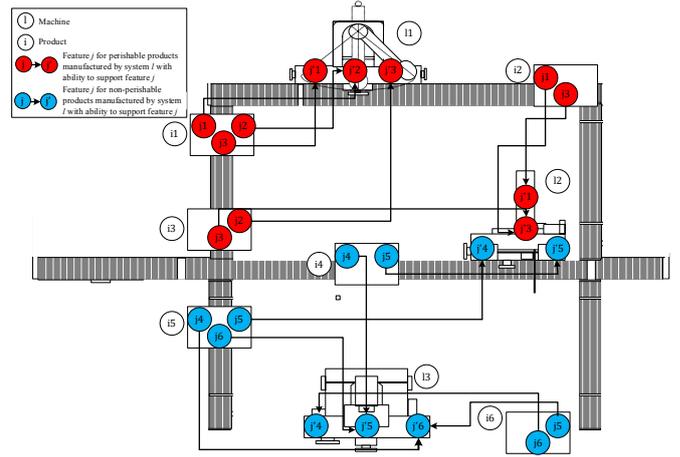

**Figure 2.** Assignment of products to production systems according to the characteristics of the products and the systems.

**Table 1.** Experimental results with the AUGMENCON2 and AUGMENCON2VIKOR algorithms

| Grid Points | AUGMENCO2 | | AUGMENCON2VIKOR | |
|---|---|---|---|---|
| | Cost | Quality | Cost | Quality |
| 1 | 15195 | 0.72 | 14665.491 | 0.735 |
| 2 | 16200 | 0.77 | 15904.654 | 0.787 |
| 3 | 17715 | 0.84 | 17336.290 | 0.859 |
| 4 | 19215 | 0.91 | 18882.404 | 0.921 |
| 5 | 19725 | 0.94 | 19269.268 | 0.956 |
| 6 | 20220 | 0.96 | 19341.654 | 0.969 |
| 7 | - | - | 19428.518 | 0.974 |
| 8 | - | - | 19515.382 | 0.978 |

*5.2. Sensitivity Analysis*

In order to demonstrate the impact of machine capacity on the assignment, machine capacity is changed. As depicted in Fig. 3, a 50% increase in machine capacity resulted an average reduction of manufacturing costs and green impacts by approximately 12% and 17%, respectively. Moreover, there was an enhancement of almost 1% in product quality.

The results of the presented algorithm, both with and without accounting green impacts are represented in Table 2. The utilization of modern machinery has a positive effect on cost and quality objective functions. For instance, in grid point 4, total cost decreased by 5429.67 units, and quality increased by 1.2%.

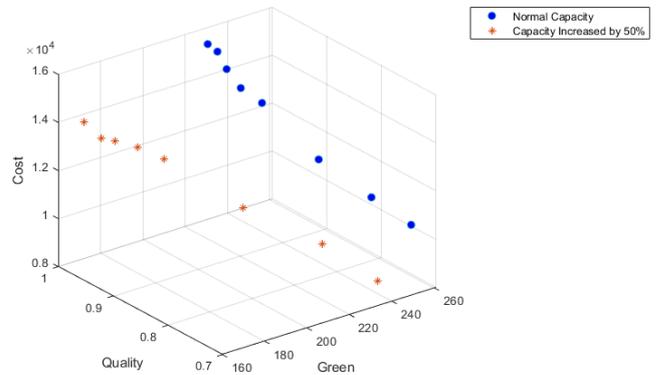

**Figure 3.** Comparative outcomes of the model under varying capacities.

Table 2. Comparative outcomes with and without environmental considerations.

| Grid Point | AUGMENCON2VIKOR without green consideration | | AUGMENCON2VIKOR with green consideration | | |
|---|---|---|---|---|---|
| | Cost | Quality | Cost | Quality | Green |
| 1 | 14665.491 | 0.735 | 10028.47 | 0.746 | 261.1 |
| 2 | 15904.654 | 0.787 | 10707.32 | 0.797 | 255.5 |
| 3 | 17336.290 | 0.859 | 11735.89 | 0.861 | 247.1 |
| 4 | 18882.404 | 0.921 | 13452.73 | 0.932 | 238.7 |
| 5 | 19269.268 | 0.956 | 13810.04 | 0.960 | 235.9 |
| 6 | 19341.654 | 0.969 | 14491.28 | 0.975 | 233.1 |
| 7 | 19428.518 | 0.974 | 15226.09 | 0.981 | 230.3 |
| 8 | 19515.382 | 0.978 | 15529.67 | 0.988 | 227.5 |

## 5. CONCLUSION

In the presented study, a novel model with three objective functions was employed to address a complex problem in sustainable manufacturing system. The products were categorized into two main categories (i.e., perishable and non-perishable products), while the machines were categorized into three main groups: dedicated machines for perishable and non-perishable goods and hybrid machines capable of performing both tasks. Each group of products had separate features and also each group of manufacturing machines had their production features. To validate the model and demonstrate its applicability, a numerical example was presented. To solve and find Pareto solution points, an AUGMECON2VIKOR algorithm was utilized for the first time in this problem. The experimental findings demonstrated the superiority of this approach over the AUGMECON2. Moreover, the inclusion of green considerations in this study resulted in notable improvements in both cost and quality with average decreases and increases of 27% and 0.8%, respectively. Investigating the problem with multiple periods and considering additional inherent complications arising from uncertainty in real sustainable manufacturing systems is an interesting and challenging direction for future work. Extending heuristic/meta-heuristic algorithms for dealing with such complicated issues is another intriguing research direction.